\documentclass[12pt]{article}
\usepackage{amsmath,amsfonts,times}
\advance\hoffset -0.5in

\newtheorem{Theorem}{Theorem}
\newtheorem{theorem}{Theorem}[section]
\newtheorem{lemma}[theorem]{Lemma}

\def\Box{{\vrule height 1.3ex width 1.0ex depth -.2ex} \hskip 1.5 truecm} 
\def\eps{\varepsilon}
\def\HH{{\mathbb H}}
\def\I{\mathcal I}

\def\RR{{\mathbb R}}
\def\SS{{\mathbb S}}
\def\XX{{\mathbb X}}
\def\ZZ{{\mathbb Z}}
\begin{document}

\author{Almut Burchard\thanks{University of Toronto, 
Toronto, Canada M5S 2E4.
{\tt almut@math.utoronto.ca}.}\quad and 
\quad Hichem Hajaiej\thanks{University of Virginia, 
Charlottesville, VA~22904-4137, USA. 
{\tt hichem.hajaiej@gmail.com}.}}

\title{~~\\[-1.5cm] Rearrangement inequalities for functionals\\
with monotone integrands}
\date{June 2004; final revision March 2006.}

\maketitle

\begin{center}
{\em Dedicated to Albert Baernstein, II
on the occasion of his 65$^{th}$ birthday.}
\end{center}

\thispagestyle{empty}
\abstract {
The 
inequalities of Hardy-Littlewood and Riesz say that certain 
integrals involving products of two or three functions increase under 
symmetric decreasing rearrangement. It is known that these inequalities  
extend to integrands of the form $F(u_1,\dots, u_m)$ where $F$ 
is supermodular; in particular, they hold
when $F$ has nonnegative mixed second derivatives
$\partial_i\partial_j F$ for all $i\ne j$.
This paper concerns the regularity assumptions on $F$ and
the equality cases.  It is shown here that 
extended Hardy-Littlewood and Riesz inequalities are valid
for supermodular integrands that are just Borel measurable.
Under some nondegeneracy conditions, all equality cases
are equivalent to radially decreasing functions
under transformations that leave the functionals
invariant (i.e., measure-preserving maps for the Hardy-Littlewood
inequality, translations for the Riesz inequality).
The proofs rely on monotone changes of variables in the spirit
of Sklar's theorem.
}

\section{Introduction} 
\label{eq:intro}


The systematic study of rearrangements
begins with the final chapter of ``Inequalities'' by Hardy, Littlewood, 
and P\'olya~\cite{HLP}. Two inequalities are discussed there 
at length, the {\em Hardy-Littlewood inequality}~(Theorems 368-370
and~378 of~\cite{HLP})
\begin{equation} \label{eq:HL}
\int_{\RR} u(x) v(x) \, dx \le \int_{\RR} u^*(x)v^*(x)\, dx\ ,
\end{equation}
and the {\em Riesz rearrangement inequality} (\cite{Riesz,Zygmund},
Theorem 370 of~\cite{HLP})
\begin{equation}
\label{eq:Riesz}
\int_{\RR}\int_{\RR} u(x)v(x')w(x-x') \, dxdx'
\le \int_{\RR}\int_{\RR} u^*(x)v^*(x')w^*(x-x')\, dxdx'\ .
\end{equation}
Here, $u$, $v$, and $w$ are nonnegative measurable functions that 
vanish at infinity, and $u^*$, $v^*$, and $w^*$ are their 
symmetric decreasing rearrangements. 

The Hardy-Littlewood inequality is a very
basic inequality that holds, with suitably defined 
rearrangements, on arbitrary measure spaces~\cite{CT}. Its main 
implication is that rearrangement decreases $L^2$-distances~\cite{Kawohl}. 
In contrast, the Riesz rearrangement inequality 
is specific to $\ZZ$ and to $\RR^n$, where it is closely
related with the Brunn-Minkowski inequality of convex geometry.
The generalization of Eq.~(\ref{eq:Riesz}) 
from $\RR$ to $\RR^n$ is due to Sobolev~\cite{Sobolev},
and the inequality is also known as the Riesz-Sobolev inequality.
For many applications, the third function in Eq.~(\ref{eq:Riesz}) 
is already radially decreasing, i.e., $w(x-x')=K(|x-x'|)$ with some 
nonnegative nonincreasing function $K$, such as the heat kernel 
or the Coulomb kernel (Theorems 371-373 and~380 of~\cite{HLP}). 
This special case of the inequality also holds on the standard 
spheres and hyperbolic spaces~\cite{BT,Beckner92}, and it still
contains the isoperimetric inequality as a limit.

It is a natural question whether these inequalities carry over to 
more general integral functionals. Under what conditions 
on $F$ do the {\em extended Hardy-Littlewood inequality}
\begin{equation}
\label{eq:HL-ext}
\int F(u_1(x),\dots, u_m(x))\, dx
\le \int F(u_1^*(x),\dots, u_m^*(x))\, dx
\end{equation}
and the {\em extended Riesz inequality}
\begin{equation}
\label{eq:Riesz-ext}
\begin{array}{l}
\displaystyle{
\hskip -.5cm
\int\cdots \int F(u_1(x_1),\dots , u_m(x_m))\, 
   \prod_{i<j} K_{ij}(d(x_i,x_j))\, dx_1\dots dx_m }\\
\displaystyle{
\qquad \le
\int\cdots \int F(u_1^*(x_1),\dots , u_m^*(x_m))\, 
   \prod_{i<j} K_{ij}(d(x_i,x_j))\, dx_1\dots dx_m }
\end{array}
\end{equation}
hold for all choices of $u_1,\dots, u_m$? In Eq.~(\ref{eq:Riesz-ext})
the $K_{ij}$ are given  nonnegative nonincreasing functions
on $\RR_+$, and $d(x,y)$ denotes the distance between $x$ and $y$. 
Eq.~(\ref{eq:HL-ext}) can be recovered from
Eq.~(\ref{eq:Riesz-ext}) by choosing 
$K_{ij}$ as a Dirac sequence and passing to the limit.  
Note that Eq.~(\ref{eq:Riesz-ext}) contains only the 
case of Eq.~(\ref{eq:Riesz}) where the third function
is a symmetric decreasing kernel. A larger class of
integral kernels $K(x_1,\dots, x_m)$
was considered in~\cite{Draghici05}.
The full generalization 
of Riesz' inequality to products of more than three functions 
was found by Brascamp-Lieb-Luttinger~\cite{BLL}; 
again, one may ask to what class of integrands 
the Brascamp-Lieb-Luttinger inequality naturally extends.  

The main condition on $F$ was identified by Lorentz~\cite{Lorentz}
as the second-order monotonicity property 
\begin{equation}
\label{eq:Lorentz}
F({\bf y } + h{\bf e}_i + k {\bf e}_j) + F({\bf y }) 
\ge F({\bf y } + h{\bf e}_i ) + F({\bf y } + k {\bf e}_j)
\qquad (i\ne j, \ h,k> 0)\ ,
\end{equation}
where ${\bf y}=(y_1,\dots, y_m)$, and ${\bf e}_i$ denotes the 
$i$-th standard basis vector in $\RR^m$.  Functions 
satisfying Eq.~(\ref{eq:Lorentz}) are called {\em supermodular} 
or {\em $2$-increasing} in Economics. A smooth function is 
supermodular, if all its mixed second partial derivatives are nonnegative.  
Eqs.~(\ref{eq:HL-ext}) and~(\ref{eq:Riesz-ext})
were proved for continuous supermodular 
integrands depending on $m=2$ functions
by Crowe-Zweibel-Rosenbloom~\cite{CZR} and Almgren-Lieb 
(Theorem 2.2 of~\cite{AL}).  For $m>2$, Eq.~(\ref{eq:HL-ext}) is 
due to Brock~\cite{Brock}, 
and Eq.~(\ref{eq:Riesz-ext}) is a recent result
of Draghici~\cite{Draghici}. The purpose of this paper is to 
dispense with the continuity assumptions on $F$ in the theorems 
of Brock and Draghici, and to characterize the equality cases 
in some relevant situations.  This continues prior work of the 
second author~[16-19].

\bigskip\noindent {\bf Acknowledgments.} \ 
We thank Friedemann Brock, Cristina Draghici, and Loren
Pitt for useful discussions, 
and especially Al Baernstein for drawing our attention to
Sklar's theorem. A.B. was partially funded by 
grants from the National Science Foundation (NSF),
the National Sciences and Engineering
Research Council of Canada (NSERC), and a University 
of Toronto Connaught award.  H.H. was supported by the Fonds National Suisse 
de la Recherche Scientifique (FNS).

\section{Statement of the results} 
\label{sec:results}
\setcounter{equation}{0}

Let $\XX$ denote either the 
Euclidean space $\RR^n$, the sphere $\SS^n$, or
the hyperbolic space $\HH^n$, equipped with the standard 
distance function $d(\cdot,\cdot)$ and the uniform volume measure $\lambda$.
Choose a distinguished point $x^*\in \XX$ to serve as
the origin or the north pole.
Consider a nonnegative measurable function $u$ on $\XX$.
When $\XX=\RR^n$ or $\HH^n$, we require $u$ to {\em vanish
at infinity} in the sense that all its positive level sets 
$\{x\in \XX:\ u(x)>t\}$ have finite measure;
when $\XX=\SS^n$ this requirement is void.
By definition, the {\em symmetric decreasing rearrangement} 
$u^*$ of $u$ is the  unique upper semicontinuous, nonincreasing function 
of $d(x,x^*)$ that is equimeasurable with~$u$.  Explicitly, if
$$
\rho(t) = \lambda\bigl(\{x\in\XX:\ u(x)>t\}\bigr)
$$
is the distribution function of $u$, and $B_r$ denotes 
the open ball of radius $r$ centered at $x^*$, then
$$
u^*(x) := \sup\bigl\{ t\ge 0:\ 
\rho(t)\ge 
 \lambda\bigl(B_{d(x,x^*)}\bigr) \bigr\} \ .
$$

\begin{Theorem} \ \label{thm:HL} 
{\bf (Extended Hardy-Littlewood inequality.)} \ 
Eq.~(\ref{eq:HL-ext}) holds for all nonnegative measurable 
functions $u_1,\dots, u_m$ that vanish at infinity
on $\XX=\RR^n,\SS^n$, or $\HH^n$, provided that the 
integrand $F$ is a supermodular Borel measurable function on 
the closed positive cone $\RR^m_+$ with $F({\bf 0})=0$,
and that its negative part satisfies
\begin{equation} \label{eq:HL-integrable}
\int_\XX F_-\bigl(u_i(x)\,{\bf e}_i\bigr)\, dx <\infty
\end{equation}
for $i=1,\dots, m$.

Suppose Eq.~(\ref{eq:HL-ext}) holds with equality, and the 
integrals are finite.  If $F$ satisfies 
Eq.~(\ref{eq:Lorentz}) with strict inequality for some $i\ne j$, 
all ${\bf y}\in\RR^m_+$ and all $h,k>0$, then
$$
\bigl(u_i(x)-u_i(x')\bigr)\bigl(u_j(x)-u_j(x')\bigr)\ge 0
$$
for almost all $x,x'\in \XX$; in particular, if $u_i=u_i^*$ is 
strictly radially decreasing, then $u_j=u_j^*$.
\end{Theorem}

\noindent The Borel measurability of $F$ and the integrability
assumption in Eq.~(\ref{eq:HL-integrable})
ensure that the integrals in Eq.~(\ref{eq:HL-ext})
are well-defined, though they may take the value $+\infty$.  

The 
left hand side of Eq.~(\ref{eq:HL-ext}) is invariant under 
volume-preserving diffeomorphisms
of $\XX$. More generally, if
$(\Omega,\mu)$ and $(\Omega',\mu')$ are measure spaces
and $\tau:\Omega\to \Omega'$ pushes $\mu$ forward to $\mu'$
in the sense that
$\mu'(A)=\mu(\tau^{-1}(A))$ for
all $\mu'$-measurable subsets $A\subset \Omega'$, then
$$
\int_\Omega  F(u_1(\omega),\dots, u_m(\omega))\, d\mu(\omega) =
\int_{\Omega'}  F(u_1\circ \tau (\omega'),\dots, u_m\circ\tau (\omega'))\, 
                   d\mu'(\omega')\ .
$$
The right hand side of Eq.~(\ref{eq:HL-ext}) can also be
expressed in an invariant form. Define the {\em nonincreasing rearrangement} 
$u^\#$ of $u$ as the unique nonincreasing upper semicontinuous 
function on $\RR_+$ that is equimeasurable with $u$,
$$
u^\#(\xi) := \sup\bigl\{ t\ge 0:\ 
\rho(t) \ge \xi\bigr\} \ .
$$
By construction, $(u\circ\tau)^\#=u^\#$ for 
any map $\tau: \Omega\to \Omega'$ that pushes $\mu$ 
forward to $\mu'$. On $\XX=\RR^n$, $\SS^n$ and $\HH^n$,
the nonincreasing rearrangement is related with 
the symmetric decreasing rearrangement by 
$u^*(x)=u^\#\bigl(\lambda\bigl(B_{d(x,x^*)}\bigr)\bigr)$.
Theorem~\ref{thm:HL} implies that
\begin{equation}
\label{eq:HL-Omega}
\int_\Omega F(u_1(\omega), \dots , u_m(\omega) )\, d\mu(\omega)
\le \int_0^{\mu(\Omega)} F(u_1^\#(\xi), \dots , u_m^\#(\xi) )\, d\xi\ 
\end{equation}
for all nonnegative measurable functions $u_1,\dots, u_m$
on $\Omega$ that vanish at infinity. 

When $\mu$ is a probability measure, Eq.~(\ref{eq:HL-Omega}) says
that the expected value of $F(Y_1,\dots, Y_m)$
is maximized among all random variables 
$Y_1,\dots, Y_m$ with given marginal distributions 
by the perfectly correlated random variables $Y_1^\#,\dots, Y_m^\#$.
The joint distribution of the maximizer is uniquely determined, if 
$Y_i$ is continuously distributed for some $i$
and Eq.~(\ref{eq:Lorentz}) is strict for all $j\ne i$.
In this formulation, the invariance under measure-preserving 
transformations is evident, since the expected value depends only on 
the joint distribution of $Y_1,\dots,Y_m$.  The assumption that 
$F$ is supermodular signifies that each of the random variables 
enhances the contribution of the others.

\begin{Theorem} \label{thm:Riesz} {\bf (Extended Riesz inequality.)} \ 
Eq.~(\ref{eq:Riesz-ext}) holds for
all nonnegative measurable functions
$u_1,\dots, u_m$ on $\XX=\RR^n,\SS^n$, or $\HH^n$ that vanish
at infinity, provided that
$F$ is a supermodular Borel measurable function 
on $\RR^m_+$ with $F({\bf 0})=0$,
each $K_{ij}$ is nonincreasing and nonnegative,
and the negative part of $F$ satisfies
\begin{equation} \label{eq:Riesz-integrable}
\int_\XX\dots \int_\XX F_-\bigl(u_\ell(x_\ell)\,{\bf e}_\ell\bigr)\ 
\prod_{i<j} K_{ij}(d(x_i,x_j))\, 
dx_1\dots dx_m <\infty
\end{equation}
for $\ell=1,\dots, m$.

Suppose Eq.~(\ref{eq:Riesz-ext}) holds with
equality. Assume additionally that the integrals are finite, 
and that $K_{ij}(t)>0$ for all $i<j$ and all $t<{\rm diam}\,\XX$.
Let $\Gamma_0$ be the graph on 
the vertex set $\{1,\dots, m\}$ which has
an edge between $i$ and $j$ whenever
$K_{ij}$ is a strictly decreasing function, and
let $i\ne j$ be from the same component 
of $\Gamma_0$.  If Eq.~(\ref{eq:Lorentz}) is strict
for all ${\bf y}\in\RR^m_+$ and all $h,k>0$,
and if $u_i$ and $u_j$ are non-constant, then 
$u_i =u_i^*\circ\tau$ and $u_j=u_j^*\circ \tau$
for some translation $\tau$ on~$\XX$. 
\end{Theorem}

\section{Related work}
\label{sec:related}
\setcounter{equation}{0} 

There are several proofs of the extended Hardy-Littlewood 
inequality in the literature. For continuous
integrands, Lorentz showed by discretization 
and elementary manipulations of the $u_i$ 
that Eq.~(\ref{eq:HL-Omega}) holds for all measurable functions 
$u_1,\dots, u_m$ on $\Omega=(0,1)$ if and only if $F$
is supermodular~\cite{Lorentz}. 
By the invariance under measure-preserving transformations,
this implies Eq.~(\ref{eq:HL-ext}), as well
as Eq.~(\ref{eq:HL-Omega}) for arbitrary finite measure 
spaces $\Omega$. However, Lorentz' paper has had little impact on 
subsequent developments.

More than thirty years later, Crowe-Zweibel-Rosenbloom proved 
Eq.~(\ref{eq:HL-ext}) for $m=2$ on $\XX=\RR^n$~\cite{CZR}.  
They expressed a given continuous supermodular function $F$ 
on $\RR^2_+$ that vanishes on the boundary as the distribution 
function of a Borel measure~$\mu_F$,
$$
F(y_1,y_2) = \mu_F\bigl( [0,y_1)\times [0,y_2)\bigr)\ .
$$
{\em layer-cake representation}

\begin{equation}
\label{eq:layercake}
\int F(u_1(x),u_2(x))\, dx
= \int_{\RR^2_+}  
\left\{ \int
\mathbf{1}_{u_1(x)> y_1} \mathbf{1}_{u_2(x)> y_2}\,
dx\right\} \, d\mu_F(y_1,y_2)\ ,
\end{equation}
which reduces Eq.~(\ref{eq:HL-ext}) to the case where $F$ is a 
product of characteristic functions~(see Theorem 1.13 in~\cite{Lieb-Loss}). 
Another reduction to products was proposed by Tahraoui~\cite{Tah}.
The regularity and boundary conditions on $F$ were
relaxed by Hajaiej-Stuart, who assumed it to be 
supermodular, of Carath\'eodory type (i.e., Borel measurable in 
the first, continuous in the second variable), and 
to satisfy some growth and integrability restrictions~\cite{HS1}. 
Equality statements for their results were obtained
by Hajaiej~\cite{H1,H2}. Using a slightly different layer-cake 
decomposition, Van Schaftingen-Willem recently 
established Eq.~(\ref{eq:HL-Omega}) for $m=2$,
under additional assumptions 
on $F$, for any equimeasurable rearrangement that preserves 
inclusions~\cite{vanS-W}. 

The drawback of the layer-cake representation is that 
for $m>2$ it requires an $m$-th order monotonicity 
condition on the integrand,
which amounts for smooth $F$ to the nonnegativity of 
all (non-repeating) mixed partial derivatives~\cite{HS2}.
Brock proved Eq.~(\ref{eq:HL-ext}) under the much weaker 
assumption that $F$ is continuous and supermodular~\cite{Brock}. 

Carlier viewed maximizing the left hand side of
Eq.~(\ref{eq:HL-Omega}) for a given right hand side
as an optimal transportation problem
where the distribution functions of $u_1,\dots, u_m$
define mass distributions $\mu_i$ on $\RR$, the joint 
distribution defines a transportation plan, and 
the functional represents the cost after multiplying by a 
minus sign~\cite{Carlier}.  He showed that the functional 
achieves its maximum (i.e., the cost is minimized) when the 
joint distribution is 
concentrated on a curve in $\RR^m$ that is nondecreasing in 
all coordinate directions, and obtained Eq.~(\ref{eq:HL-Omega}) 
as a corollary.  His proof takes advantage of the dual problem 
of minimizing 
$$
\sum_{i=1}^m \int_\RR f_i(y)\,  d\mu_i(y) 
$$
over $f_1,\dots, f_m$, subject to the constraint that 
$\sum f_i(y_i)\ge F(y_1,\dots, y_m)$
for all $y_1,\dots, y_m$.

Theorem~\ref{thm:HL} can be applied to 
some integrands that depend explicitly on the radial 
variable~\cite{Lorentz,HS1,Carlier}.
If $G$ is a function on $\RR_+\times \RR^m_+$
such that
$F(y_0,\dots,y_m):= G(y_0^{-1}, y_1,\dots,y_m)$ 
satisfies the assumptions of Theorem~\ref{thm:HL}, then
\begin{equation} \label{eq:HL-x}
\int_{\RR^n} G(|x|,u_1(x),\dots, u_m(x))\, dx
\le \int_{\RR^n} G(|x|,u_1^*(x),\dots, u_m^*(x))\, dx\ .
\end{equation}
Hajaiej-Stuart studied this
inequality in connection with 
the following problem in nonlinear optics~\cite{HS1,HS2}.
The profiles of stable electromagnetic waves traveling along 
a planar wave\-guide are given 
by the ground states of the energy functional
$$
{\cal E}(u) = \frac{1}{2} \int_\RR |u'|^2\, dx - \int_\RR G(|x|,u)\, dx\ 
$$
under the constraint $||u||_2=c$.  Here, $x$ is the position relative
to the optical axis, $G$ is determined by the index of refraction, 
and $c>0$ is a parameter related to the wave speed~\cite{Stuart93}.  
If the index of refraction of the optical 
media decreases with $|x|$, then $F(r,y)=G(r^{-1},y)$ satisfies 
the assumptions of Theorem~\ref{thm:HL}. Then the first integral shrinks  
under symmetric decreasing rearrangement by the P\'olya-Szeg\H o 
inequality, the second integral grows by Eq.~(\ref{eq:HL-x}), 
and the $L^2$-constraint is conserved. Thus, one may rearrange 
any minimizing sequence to obtain a minimizing sequence of symmetric 
decreasing functions.  This is a crucial step in 
the construction of ground states --- if $G$ violates the
monotonicity conditions, then a ground state need not exist~\cite{HS3}.  
Hajaiej-Stuart worried about restrictive regularity 
assumptions, because $G$ may 
jump 
at interfaces between layers of different media.

The Riesz inequality in Eq.~(\ref{eq:Riesz-ext}) 
is non-trivial even when $F$ is just a product of
two functions.  Ahlfors introduced
two-point rearrangements to treat this case 
on $\XX=\SS^1$~\cite{Ahlfors},
Baernstein-Taylor proved the corresponding result on $\SS^n$~\cite{BT},
and Beckner noted that the proof remains valid on 
$\HH^n$ and $\RR^n$~\cite{Beckner92}.  When $F$ is a product of $m>2$ 
functions, Eq.~(\ref{eq:Riesz-ext}) has applications to spectral 
invariants of heat kernels via the Trotter product 
formula~\cite{Luttinger}. This case
was settled by Friedberg-Luttinger~\cite{FL},
Burchard-Schmuckenschl\"ager~\cite{BS}, and by Morpurgo,
who proved Eq.~(\ref{eq:Riesz-ext})
more generally for integrands of the form
\begin{equation}
\label{eq:Morpurgo}
F(y_1,\dots, y_m)= \Phi\Bigl(\sum_{i=1}^m y_i\Bigr)
\end{equation}
with $\Phi$ convex (Theorem 3.13 of~\cite{Mor}). 
In the above situations, equality cases 
have been determined~\cite{Lieb-Choq,Beckner93,BS,Mor}.
Almgren-Lieb used the technique of Crowe-Zweibel-Rosen\-bloom
to prove Eq.~(\ref{eq:Riesz-ext}) for $m=2$~\cite{AL}.
The special case where $F(u,v)=\Phi(|u-v|)$ for some 
convex function $\Phi$ was identified
by Baernstein  as a `master inequality' from which many classical 
geometric inequalities can be derived quickly~\cite{B-uni}.  
Eq.~(\ref{eq:Riesz-ext}) for continuous supermodular integrands 
with $m>2$ is due to Draghici~\cite{Draghici}.

\section{\bf Outline of the arguments}
\label{sec:outline}
\setcounter{equation}{0} 
In their proofs of Eqs.~(\ref{eq:HL-ext}) and~(\ref{eq:Riesz-ext}), 
Brock and Draghici showed that the left hand sides 
increase under two-point rearrangements if $F$ is any supermodular 
Borel integrand~\cite{Brock,Draghici}.  Then they 
approximated the symmetric decreasing rearrangement with sequences 
of repeated two-point rearrangements. 
Baernstein-Taylor had established that such sequences
can be made to converge to the symmetric decreasing 
rearrangement in a space of continuous functions 
~\cite{BT}, and Brock-Solynin had proved this 
convergence in $L^p$-spaces~\cite{Brock-Solynin}.
To pass to the desired limits, Brock and Draghici assumed that 
$F$ is continuous and satisfies some boundary and growth conditions. 

No new proofs of these inequalities will be given here.  Rather, 
we reduce general supermodular integrands to the known cases
of integrands that are also bounded and continuous. 
This reduction needs more care than the usual density arguments, 
because pointwise $a.e.$ convergence of a sequence of integrands $F_k$ does 
not guarantee pointwise $a.e.$ convergence of the compositions
$F_k(u_1,\dots, u_m)$. Approximation within a 
class of functions with specified positivity or 
monotonicity properties can be subtle; for instance, 
nonnegative functions of $m$ variables cannot always 
be approximated by positive linear combinations of products of
nonnegative functions of the individual variables
(contrary to Theorem 2.1 and Lemma~4.1 of~\cite{Tah}). 

In Section~\ref{sec:super}, we 
prove a variant of Sklar's theorem~\cite{Sklar59} which
factorizes a given supermodular function
on $\RR^m_+$ as the composition of a Lipschitz continuous 
supermodular function on $\RR^m_+$ with $m$ 
monotone functions on $\RR_+$, and a cutoff lemma 
that replaces a given supermodular function by a bounded
supermodular function. Section~\ref{sec:2pt} is dedicated to 
the two-point versions of Theorems~\ref{thm:HL} and~\ref{thm:Riesz}.
Here, we review the proofs of the two-point rearrangement 
inequalities of Lorentz~\cite{Lorentz}, Brock~\cite{Brock}, and 
Draghici~\cite{Draghici} and find their equality cases.
The main theorems are proved in Section~\ref{sec:proofs}
by combining the results from Sections~\ref{sec:super} 
and~\ref{sec:2pt}. Adapting Beckner's argument from~\cite{Beckner93}, 
we note that the inequalities in Eq.~(\ref{eq:HL-ext}) and 
Eq.~(\ref{eq:Riesz-ext}) are strict unless $u_1,\dots, u_m$ produce 
equality in {\em all} of the corresponding two-point inequalities,
and then apply the results from Section~\ref{sec:2pt}.
In the final Section~\ref{sec:conc}, we briefly discuss
extensions for the Brascamp-Lieb-Luttinger and 
related inequalities.

\section{Monotone functions} 
\label{sec:super}
\setcounter{equation}{0}

In this section, we provide  two technical results about 
functions with higher-order monotonicity properties.  
We begin with an auxiliary lemma for functions of a 
single variable.

\begin{lemma} \label{lem:factorize-1} {(\bf Monotone change of variable.)} \ 
Let $\phi$ be a nondecreasing
real-valued function defined on an interval $I$.
Then, for every function $f$ on $I$ satisfying 
\begin{equation}
\label{eq:factorize-1b}
|f(z)-f(y)| \le C(\phi(z)-\phi(y))
\end{equation}
for all points $y<z\in I$ with some constant $C$,
there exists a Lipschitz continuous function
$\tilde f: \RR\to [\inf f,\sup f]$ 
such that $f=\tilde f\circ \phi$. Furthermore, if
$f$ is nondecreasing, then $\tilde f$ is nondecreasing.
\end{lemma}

\noindent{\sc Proof.} 
If $t=\phi(y)$ we set $\tilde f(t) := f(y)$.
For $s<t$ with $s=\phi(y)$, $t=\phi(z)$,
Eq.~(\ref{eq:factorize-1b}) implies that
\begin{equation}
\label{eq:factorize-1a}
|\tilde f(t)-\tilde f(s)| = |f(z)-f(y)|\le C(\phi(z)-\phi(y))
= C(t-s)\ . 
\end{equation}
Since $\tilde f$ is uniformly continuous on the image
of $\phi$, it has a unique continuous extension
to the closure of the image.  The complement 
consists of a countable number of open disjoint bounded intervals, 
each representing
a jump of $\phi$, and possibly one or two unbounded
intervals. On each of the bounded intervals, we interpolate 
$\tilde f$ linearly between the values that have already been 
assigned at the endpoints.
If $\phi$ is bounded either above or below, we extrapolate 
$\tilde f$ to $t>\sup \phi$ and $t<\inf \phi$ by constants.

By construction,  $f=\tilde f\circ \phi$ and 
$\tilde f(\RR)=[\inf f,\sup f]$. The continuous extension 
and the linear interpolation preserve the modulus 
of continuity of $\tilde f$, and hence, by Eq.~(\ref{eq:factorize-1a}), 
\begin{equation}
|\tilde f(t)-\tilde f(s)|\le C|t-s|
\end{equation}
for all $s,t\in\RR$. If $f$ is nondecreasing, then $\tilde f$ is 
nondecreasing on the image of $\phi$ by definition, and on the 
complement by continuous extension and linear interpolation.  
\hfill $\Box$\quad

\bigskip Lemma~\ref{lem:factorize-1} is related to the
elementary fact that a continuous random variable 
can be made uniform by a monotone change of variables.
More generally, if $\phi$ is nondecreasing and right 
continuous, and its generalized inverse is
defined by $\psi(t) = \inf\{ y:\ \phi(y) \ge  t \}$,
then the cumulative distribution functions
of two random variables that are related by $Y=\psi(\tilde Y)$ satisfy
$$
F(y) = P(Y\le y) = P\bigl( \tilde Y\le \phi(y)\bigr) 
= \tilde F\bigl(\phi(y)\bigr)\ ,
$$
i.e., $F=\tilde F\circ \phi$.  Choosing $\phi=F$ results in 
a uniform distribution for $\tilde Y$.

The corresponding result for $m\ge 2$ random variables is
known as {\em Sklar's theorem}~\cite{Sklar59}. The theorem asserts
that a collection of random variables $Y_1,\dots, Y_m$
with a given joint distribution function
$F$ can be replaced by random variables $\tilde Y_1,\dots, \tilde Y_m$
whose marginals $\tilde Y_i$ are uniformly distributed
on $[0,1]$, and whose joint distribution function $\tilde F$ 
is continuous. The next lemma contains Sklar's 
theorem for supermodular functions. Since the lemma 
follows from the arguments outlined in~\cite{Sklar73} rather than 
from the statement of the theorem, we include its proof
for the convenience of the reader.

We first introduce some notation. Let $F$ be a real-valued function on 
the closed positive cone $\RR^m_+$.  For $i=1,\dots, m$ and 
$h\ge 0$, consider the finite difference operators
$$
\Delta_i F({\bf y};h) :=  F({\bf y}+h{\bf e}_i)-F({\bf y})\ .
$$
The operators commute, and higher order difference operators 
are defined recursively by
$$
\Delta_{i_1\dots i_\ell} F({\bf y}; h_1,\dots, h_\ell)
:= \Delta_{i_1\dots i_{\ell-1}} 
  \Delta_{i_\ell} F(({\bf y}; h_\ell); h_1,\dots, h_{\ell-1})\ .
$$
If $F$ is $\ell$ times continuously differentiable, then
$$
\Delta_{i_1\dots i_\ell} F({\bf y}; h_1,\dots, h_\ell)
= \int_0^{h_1}\cdots\int_0^{h_\ell}
\partial_{i_1}\dots\partial_{i_\ell} F\Bigl({\bf y} + \sum_{i=1}^\ell
t_i{\bf e}_i\Bigr)\, dt_1\dots dt_\ell\ .
$$ 
A function $F$ is nondecreasing in each variable if
$\Delta_iF\ge 0$ for $i=1,\dots, m$; it is supermodular, 
if $\Delta_{ij} F\ge 0$ for all $i\ne j$.  
The joint distribution function of $m$ random variables satisfies 
$\Delta_{i_1}\dots\Delta_{i_\ell}F\ge 0$ for any choice of
distinct indices $i_1,\dots, i_\ell$.

\begin{lemma} {\bf (Sklar's theorem.)} \  \label{lem:Sklar} 
Assume that $F$ is boun\-ded, nondecreasing
in each variable, and supermodular on $\RR^m_+$.  Then there 
exist bounded nondecreasing functions $\phi_1,\dots, \phi_m$
on $\RR_+$ with $\phi_i(0)=0$ and a Lipschitz continuous
function $\tilde F$ on $\RR^m_+$  such that
$$
F(y_1,\dots, y_m)= \tilde F(\phi_1(y_1), \dots,\phi_m(y_m))\ .
$$
Furthermore, $\tilde F$ is bounded, nondecreasing in each
variable, and supermodular.
If, in addition, $\Delta_{i_1\dots i_\ell}F\ge 0$ 
on $\RR^m_+\times \RR_+^\ell$ for 
some distinct indices $i_1,\dots, i_\ell$,
then  $\Delta_{i_1\dots i_\ell}\tilde F\ge 0$.
\end{lemma}

\noindent{\sc Proof.} 
Set
$$
\phi_i(y) = \lim_{y_j\to\infty, j\ne i} 
\left\{ F(y_1,\dots, y_m)\Big\vert_{y_i=y}  - 
F(y_1,\dots, y_m)\Big\vert_{y_i=0}\right\} \ .
$$
These functions are nonnegative and bounded by $\sup F-\inf F$. Since $F$ 
is nondecreasing in each variable, they 
are nonnegative, and since $F$ is supermodular, they are
nondecreasing and satisfy
\begin{equation} \label{eq:assumption-factorize}
F({\bf y}+ h{\bf e}_i)- F({\bf y})\le \phi_i(y_i+h)-\phi_i(y_i)\ 
\end{equation}
for all ${\bf y}=(y_1,\dots,y_m)\in \RR^m_+$ and all $h>0$.

We construct $\tilde F$ by changing 
one variable at a time. For the first variable,
we write ${\bf y}=(y,\hat {\bf y})$ where $y\in\RR_+$ and 
$\hat {\bf y}\in \RR^{m-1}_+$. 
By Eq.~(\ref{eq:assumption-factorize}), for each 
$\hat {\bf y}\in \RR^{m-1}_+$, the function
$f(y)= F(y,\hat {\bf y})$ satisfies Eq.~(\ref{eq:factorize-1b})
with $C=1$ and $\phi=\phi_1$.
By Lemma~\ref{lem:factorize-1}, 
there exists a function $F_1$ satisfying 
$$
F(y,\hat {\bf y}) = F_1(\phi_1(y),\hat {\bf y})
$$
for all $(y,\hat {\bf y})\in\RR^m_+$. Furthermore, $F_1$ is 
Lipschitz continuous in the first variable,
$$
| F_1(t,\hat {\bf y})-F_1(s,\hat {\bf y})|\le |t-s|\ .
$$
We claim that $F_1$ satisfies
Eq.~(\ref{eq:assumption-factorize}) for all $j>1$ 
with the same function $\phi_j$ as $F$.
To see this, note that for each $h>0$ and every $\hat {\bf y}$, 
$$
f(y) = \Delta_jF(y,\hat {\bf y};h)
$$
satisfies the assumptions of Lemma~\ref{lem:factorize-1}
with $C=2$ and $\phi=\phi_1$. A moment's consideration shows that 
$$
\tilde f(t) =  \Delta_j F_1(t,\hat {\bf y}; h)
$$
and the claim follows since $\sup \tilde f =\sup f\le 
\phi_j(y_j+h)-\phi_j(y_j)$ by Lemma~\ref{lem:factorize-1}.

We next verify that $F_1$ has the same monotonicity properties 
as $F$.  Suppose that $\Delta_{i_1\dots i_\ell}F\ge 0$
for some set of $\ell\ge 1$ distinct indices 
$i_1,\dots, i_\ell$.  If $1\not\in \{i_1,\dots, i_\ell\}$, 
we apply Lemma~\ref{lem:factorize-1} to 
$f(y)=\Delta_{i_1\dots i_\ell}F(y,\hat {\bf y}; h_1,\dots, h_\ell)$,
which satisfies Eq.~(\ref{eq:factorize-1b}) with $C=2^\ell$ and $\phi=\phi_1$
for all $\hat {\bf y}\in\RR^{m-1}$ and all $h_1,\dots, h_\ell\ge 0$.
It follows that
$\tilde f(t)=\Delta_{i_1\dots i_\ell} F_1(t,\hat {\bf y}; h_1,\dots, h_\ell)
\ge 0$.  On the other hand, if $i_1=1$, we 
apply Lemma~\ref{lem:factorize-1} to
$f(y)=\Delta_{i_2,\dots, i_\ell} F(y,\hat {\bf y}; h_2,\dots, h_\ell)$.
Since $f(y)$ is nondecreasing by assumption, $\tilde f(t) 
=\Delta_{i_2,\dots, i_\ell} F_1(t,\hat {\bf y}; h_2,\dots, h_\ell)$ 
is again nondecreasing, and we conclude that 
$\Delta_{i_1\dots i_\ell} F_1\ge 0$ also in this case.

Iterating the change of variables for $i=2,\dots, m$
gives functions $F_i$ satisfying
$$
F_{i-1}(t_1,\dots, t_{i-1}, y_i, \dots, y_m)
= F_i(t_1,\dots,t_{i-1},\phi_i(y_i),y_{i+1},\dots y_m)\ ,
$$
as well as
\begin{equation} \label{eq:factorize-Lip-i}
0\le \Delta_j F_i(t_1,\dots,t_i,y_{i+1},\dots y_m ; h)
\le \left\{\begin{array}{ll} 
h\ , \quad & j\le i\\
\phi_j(y_j+h)-\phi_j(y_j)\ , \quad & j> i\ .
\end{array}\right.
\end{equation}

Finally, we set $\tilde F=F_m$.
It follows from Eq.~(\ref{eq:factorize-Lip-i})
that $\tilde F$ satisfies the Lipschitz condition
$|\tilde F({\bf z}) - \tilde F({\bf y})| \le \sum |z_i-y_i|
\le \sqrt{m}\, |{\bf z}-{\bf y}|$ for all 
${\bf y}, {\bf z}\in\RR^m_+$. \hfill\Box\quad

\bigskip\noindent
The distribution function of a Borel measure on $\RR^m_+$ can be
conveniently approximated from below by restricting 
the measure to a large cube $[0,L)^m$. The next 
lemma constructs the corresponding approximation 
for functions with weaker monotonicity properties.

\begin{lemma} \label{lem:cutoff} {\bf (Cutoff.)} \ 
Given a real-valued function $F$ in $\RR^m_+$, set
$$
F^L(y_1,\dots, y_m) := F(\min\{y_1,L\},\dots, \min\{y_m,L\})\ .
$$
If $F$ is nondecreasing in each variable, then 
$F^L\le F$.  If $\Delta_{i_1\dots i_\ell}F\ge 0$  
on $\RR^m_+\times \RR_+^\ell$ 
for some distinct indices $i_1,\dots, i_\ell$, then
$\Delta_{i_1,\dots, i_\ell}F^L\ge 0$.
In particular, if $F$ is supermodular, so is $F^L$.
If $F$ has the property that $\Delta_{i_1\dots i_\ell}F\ge 0$ 
on $\RR^m_+\times\RR_+^\ell$ for every set of distinct indices 
$i_1,\dots,i_\ell$, then $F-F^L$ also has this property.
\end{lemma}

\noindent{\sc Proof.}  As in the proof of Lemma~\ref{lem:Sklar},
we modify the variables one at a time.
The function
$
F^{1,L} (y,\hat {\bf y}) := F\bigl (\min\{y,L\}, \hat {\bf y}\bigr)$
has the same monotonicity properties as $F$
because $\min\{y,L\}$ is nondecreasing in~$y$.

If $\Delta_{i_1\dots i_\ell}F\ge 0$ for
all collections of distinct indices $i_1,\dots, i_\ell$, we write
$$
F(y,\hat {\bf y})-F^{1,L}(y,\hat {\bf y}) 
= \Delta_1 F\bigl(y,\hat {\bf y};[y-L]_+\bigr)\ ,
$$
and it follows that
$\Delta_{i_1\dots i_\ell} (F-F^{1,L})\ge 0$
whenever $1\not\in \{i_1,\dots,i_\ell\}$. 
For $i_1=1$, we write
$$
\Delta_1( F(y,\hat {\bf y}; h)
-F^{1,L} (y,\hat {\bf y}; h)
= \Delta_1F\bigl(\max\{y,L\},\hat {\bf y};[h-[L-y]_+]_+\bigr)\ ,
$$
and conclude that $\Delta_{i_1\dots i_\ell}(F-F^{1,L})\ge 0$
also in this case.  

Repeating the construction for  the variables $y_2,\dots, y_m$ 
gives the claims.
\hfill $\Box$\quad

\section{Two-point rearrangements} 
\label{sec:2pt}
\setcounter{equation}{0}

Let $\XX$ be $\RR^n$, $\SS^n$, or $\HH^n$. 
A {\em reflection} on $\XX$ is an isometry
characterized by the properties that 
\ {\em(i)}~~$\sigma^2x=x$ for all $x\in\XX$; 
\ {\em(ii)}~~the fixed point set $H_0$ of $\sigma$ separates $M$ into
two half-spaces $H_+$ and $H_-$ that are interchanged by $\sigma$; 
\ and {\em (iii)}~~$d(x,x')<d(x,\sigma x')$ for all $x,x'\in H_+$.
We call $H_+$ and $H_-$ the {\em positive} and {\em negative half-spaces}
associated with $\sigma$.  By convention, we always
choose $H_+$ to contain the distinguished
point $x^*$ of $\XX$ in its closure.
The {\em two-point rearrangement}, or {\em polarization} 
of a real-valued function $u$ with respect to a reflection 
$\sigma$ is defined by 
$$
u^\sigma(x) = \left\{ \begin{array}{ll}
\max\{u(x),u(\sigma x)\}\ , \quad & x\in H_+\cup H_0\\
\min\{u(x),u(\sigma x)\}\ , \quad & x\in H_-\; .
\end{array}\right.
$$
This definition makes sense, and
the two-point versions of Eqs.~(\ref{eq:HL-ext}) 
and~(\ref{eq:Riesz-ext}) hold for any space with a reflection 
symmetry.

On $\XX=\RR^n$, $\SS^n$, and $\HH^n$, any pair of points 
is connected by a unique reflection.
The space of reflections forms an $n$-dimensional 
submanifold of the $n(n+1)/2$-dimensional space 
of isometries, and thus has a natural uniform metric.
If $u$ is measurable, both the composition
$u\circ\sigma$ and the rearrangement $u^\sigma$ depend
continuously on $\sigma$ in the sense that
$\sigma_k\to \sigma$ implies that
$u\circ\sigma_k\to u\circ\sigma$
and $u^{\sigma_k}\to u^\sigma$ in measure.

Two-point rearrangements are particularly well-suited for identifying
symmetric decreasing functions, because
\begin{equation}
\label{eq:2pt-identify}
u = u^* \quad \Longleftrightarrow \quad u=u^\sigma\ \mbox{for all
$\sigma$}\ .
\end{equation}
Functions that are radially decreasing about some 
point are characterized by
\begin{equation}
\label{eq:2pt-identify-tau}
u = u^*\circ \tau \ \mbox{for some translation $\tau$}
\quad \Longleftrightarrow \quad  \mbox{for all $\sigma$, 
either $u=u^\sigma$ or $u=u^\sigma\circ\sigma$}
\end{equation} 
(see Lemma 2.8 of~\cite{BS}).

Integral inequalities for two-point rearrangements typically
reduce to elementary combinatorial inequalities for the 
integrands.  The following lemma supplies the
elementary inequality for the Hardy-Littlewood and Riesz functionals.

\begin{lemma}\label{lem:Lorentz} {\bf (Lorentz two-point 
inequality.)} \ 
A real-valued function $F$ on $\RR^m_+$ is supermodular,
if and only if for every pair of points ${\bf z}, {\bf w} \in\RR^m_+$.
\begin{equation}
\label{eq:algebra-claim}
\begin{array}{lcl}
F(z_1,\dots , z_m) + F(w_1,\dots , w_m) &\le &
F(\max \{z_1,w_1\},\dots,\max\{z_m,w_m\})\\
&&  + 
F(\min \{z_1,w_1\},\dots,\min\{z_m,w_m\}) \,.
\end{array}
\end{equation}
If $\Delta_{ij}F>0$ for some $i\ne j$ 
then Eq.~(\ref{eq:algebra-claim}) is strict 
unless $(z_i-w_i)(z_j-w_j)\ge 0$.
\end{lemma}

\noindent{\sc Proof.}
Given ${\bf z}, {\bf w}\in\RR^m_+$, define 
${\bf y}, {\bf h}\in \RR^m_+$ by $y_i = \min\{z_i,w_i\}$ and 
$h_i=|z_i-w_i|$  for $i=1,\dots m$.
If $I\subset \{1,\dots,m\}$, we use the notation 
${\bf h}_I = \sum_{i\in I} h_i{\bf e}_i$.
Subtracting the left hand side of Eq.~(\ref{eq:algebra-claim})
from the right hand side  results in the equivalent 
statement
\begin{equation} \label{eq:algebra-proof}
\Delta_{IJ}F({\bf y}; {\bf h}_I,{\bf h}_J):=
F({\bf y}+ {\bf h}_{I\cup J})
- F({\bf y}+ {\bf h}_I) - F({\bf y}+ {\bf h}_J)
+ F({\bf y})\ge 0\ ,
\end{equation}
where $I=\{i:\ z_i<w_i\}$, and $J=\{i: z_i>w_i\}$.
If either $I$ or $J$ is empty, Eq.~(\ref{eq:algebra-proof})
is trivially satisfied. If $I$ and $J$ each have exactly 
one element, Eq.~(\ref{eq:algebra-proof})
is equivalent to Eq.~(\ref{eq:Lorentz}).
If one of the sets, say $I$, has several elements, then decomposing 
it into disjoint subsets as $I=I'\cup I''$ gives
$$
\Delta_{IJ}F({\bf y}; {\bf h}_I, {\bf h}_J)
=\Delta_{I'J}F ({\bf y} + {\bf h}_{I''}, {\bf h}_{I'},{\bf h}_J)
+ \Delta_{I''J}F ({\bf y} , {\bf h}_{I'},{\bf h}_J)\ ,
$$
and Eq.~(\ref{eq:algebra-proof}) follows by recursion.
The same recursion implies that if
$\Delta_{ij}F>0$ and $z_i-w_i$ and $z_j-w_j$ have opposite signs, then
the inequality in Eq.~(\ref{eq:algebra-proof})
is strict whenever  $I$ contains~$i$, $J$ contains~$j$, 
and $h_i,h_j>0$. \hfill $\Box$\quad

\bigskip\noindent 
Brock proved that the left hand side of Eq.~(\ref{eq:HL-ext})
increases under two-point rearrangement~\cite{Brock}:

\begin{lemma} \label{lem:HL-2pt} 
{\bf (Hardy-Littlewood two-point inequality.)} \ 
Let $F$ be a 
supermodular Borel measurable function
on $\RR^m_+$, and let $u_1,\dots,u_m$ be nonnegative measurable functions
on $\XX$ satisfying the integrability
condition in Eq.~(\ref{eq:HL-integrable}).
Then, for any reflection $\sigma$ on $\XX$, 
\begin{equation} \label{eq:HL-2pt}
\int_{\XX} F\bigl(u_1(x),\dots , u_m(x)\bigr)\, dx
\le \int_{\XX} F\bigl(u_1^\sigma(x),\dots , u_m^\sigma(x)\bigr)\, dx\ .
\end{equation}

Assume furthermore that $\Delta_{ij}F>0$ on $\RR^m_+\times (0,\infty)^2$
for some $i\ne j$.  If Eq.~(\ref{eq:HL-2pt}) holds with equality
and the integrals are finite,  then 
$$
\bigl(u_i(x)-u_i(\sigma x)\bigr)\bigl(u_j(x)-u_j(\sigma x)\bigr)\ge 0
\quad a.e.\,.
$$
In particular, if $u_i=u_i^*$ is strictly radially 
decreasing and $\sigma(x^*)\ne x^*$, then $u_j=u_j^\sigma$.
\end{lemma}

\noindent{\sc Proof.} \ {\em The inequality~\cite{Brock}:}
The left hand side of Eq.~(\ref{eq:HL-2pt})
can be written as an integral over the positive half-space,
$$
\I(u_1,\dots, u_m) := \int_{H_+}
F\bigl(u_1(x),\dots u_m(x)\bigr) + 
F\bigl(u_1(\sigma x),\dots u_m(\sigma x)\bigr)\, dx\ .
$$
By Lemma~\ref{lem:Lorentz}, 
with $z_i=u_i(x)$ and $w_i=u_i(\sigma x)$,  
the integrand satisfies
\begin{equation} \label{eq:HL-2pt-integrand}
\begin{array}{l}
\hskip -2cm
F\bigl(u_1(x),\dots u_m(x)\bigr) + 
F\bigl(u_1(\sigma x),\dots u_m(\sigma x)\bigr)\\[0.1cm]
\le 
F\bigl(u_1^\sigma(x),\dots u_m^\sigma(x)\bigr) + 
F\bigl(u_1^\sigma(\sigma x),\dots u_m^\sigma(\sigma x)\bigr)
\end{array}
\end{equation}
for all $x\in H_+$.  Integrating over $H_+$ yields Eq.~(\ref{eq:HL-2pt}).

\medskip{\em Equality statement:} 
Assume that $\I(u_1,\dots , u_m)= \I(u_1^\sigma,\dots , u_m^\sigma)$
is finite.  Then Eq.~(\ref{eq:HL-2pt-integrand})
must hold with equality almost everywhere on $H_+$. 
If $\Delta_{ij}F>0$ on $\RR^m_+\times (0,\infty)^2$, 
then Lemma~\ref{lem:Lorentz} implies that 
$u_i(x)-u_i(\sigma x)$ and $u_j(x)-u_j(\sigma x)$
cannot have opposite signs except on a set of
zero measure.  If moreover $u_i=u_i^*$ is strictly radially decreasing
and $\sigma x^*\ne x^*$, then
$u_i(x)>u_i(\sigma x)$ for $a.e.\ x\in H_+$,
and Lemma~\ref{lem:Lorentz} implies that $u_j(x)\ge u_j(\sigma x)$ for
$a.e.\  x\in H_+$.
\hfill $\Box$\quad

\bigskip \noindent Brock completed the proof of
Eq.~(\ref{eq:HL-ext}) by approximating
the symmetric decreasing rearrangement with a sequence of two-point 
rearrangements \`a la Baernstein-Taylor~\cite{BT}. We 
sketch his argument in the simplest case where $F$ is 
a continuous supermodular function that vanishes on the boundary of 
the positive cone $\RR^m_+$, and $u_1,\dots, u_m$ are bounded 
and compactly supported.

By Theorem 6.1 of~\cite{Brock-Solynin} there exists a sequence 
of reflections $\{\sigma_k\}_{k\ge 1}$ such that 
\begin{equation}
\label{eq:2pt-converge}
u_i^{\sigma_1,\dots, \sigma_k} \to u_i^*\quad \mbox{in measure}
\ (k\to\infty)
\end{equation}
for $i=1,\dots, m$. By Lemma~\ref{lem:HL-2pt}, the 
functional increases monotonically along such 
a sequence.  If $B$ is a ball centered at $x^*$ that contains
the supports of $u_1,\dots, u_m$, then 
the rearranged functions $u_i^{\sigma_1,\dots, \sigma_k}$ are 
also supported on $B$, and dominated convergence yields
\begin{equation}
\label{eq:dominated}
\I(u_1,\dots,u_m) 
\le \I (u_1^{\sigma_1,\dots, \sigma_k},\dots , u_m^{\sigma_1,\dots, 
\sigma_k}) \to \I (u_1^*,\dots, u_m^*)\quad (k\to\infty) \ .
\end{equation}

\bigskip  
The corresponding results for Eq.~(\ref{eq:Riesz-ext}) are due to 
Draghici~\cite{Draghici}. The two-point inequality
is not an immediate consequence of Lemma~\ref{lem:Lorentz}, but 
requires an additional combinatorial argument. This argument
was used previously by Morpurgo~\cite{Mor}, 
and a simpler version appears in~\cite{BS}.

\begin{lemma} \label{lem:Riesz-2pt}
{\bf (Riesz two-point inequality.)} \ 
Assume that $F$ is a supermodular Borel measurable function on $\RR^m_+$.
For each pair of indices
$1\le i<j\le m$, let $K_{ij}$ be a  nonincreasing
function on $\RR_+$, and let
$u_1,\dots, u_m$ be nonnegative measurable functions on $\XX$
satisfying the integrability condition
in Eq.~(\ref{eq:Riesz-integrable}).
Then, for any reflection $\sigma$,
\begin{equation} \label{eq:Riesz-2pt}
\begin{array} {l}
\displaystyle{
\hskip -0.5cm\int_{\XX}\cdots \int_{\XX} F\bigl(u_1(x_1),\dots , u_m(x_m)\bigr)
\prod_{i<j} K_{ij}\bigl(d(x_i,x_j)\bigr)\, dx_1\dots dx_m}\\
\displaystyle{\quad
\le \int_{\XX} \dots \int_{\XX} 
   F\bigl(u_1^\sigma(x_1),\dots , u_m^\sigma(x_m)\bigr)
\prod_{i<j}K_{ij}\bigl(d(x_i,x_j)\bigr)\, 
dx_1\dots dx_m\ .}
\end{array}
\end{equation}

Assume additionally that that $K_{ij}(t)>0$ for all 
$i<j$ and all $t<{\rm diam}\,\XX$.  Let $\Gamma_0$ be the graph 
on $\{1,\dots,m\}$ with an edge between $i$ and $j$ whenever
$K_{ij}$ is strictly decreasing.
If $\Delta_{ij}F>0$ for some  $i\ne j$ lying in the same connected
component of $\Gamma_0$,  and that
$u_i$ and $u_j$ are not symmetric under $\sigma$.
If the integrals in Eq.~(\ref{eq:Riesz-2pt}) have the same
finite value, then either $u_i=u_i^\sigma$ 
and $u_j=u_j^\sigma$, or $u_i= u_i^\sigma\circ\sigma$ 
and $u_j= u_j^\sigma\circ\sigma$.
\end{lemma}

\noindent{\sc Proof.} \ {\em The inequality~\cite{Draghici}:} 
The left hand side of Eq.~(\ref{eq:Riesz-2pt})
can be written as an $m$-fold integral over the positive half-space 
\begin{eqnarray}
\nonumber
\I(u_1,\dots, u_m) &:=&
\int_{H_+}\dots\int_{H_+}
\sum_{\eps_i\in\{0,1\}, i=1,\dots,m} 
\biggl\{
F\bigl(u_1(\sigma^{\eps_1}x_1),\dots, u_m(\sigma^{\eps_m}x_m)\bigr)\, \times
\\
\label{eq:Riesz-2pt-proof-1}
&& \hskip 1.5cm \times \ 
\prod_{i<j} K_{ij}\bigl(d(\sigma^{\eps_i}x_i,\sigma^{\eps_j}x_j)\bigr)\biggr\}
\, dx_1 \dots dx_m\ .
\end{eqnarray}
Fix $x_1,\dots, x_m\in H_+$. For each pair of indices $i<j$, set
$a_{ij}= K_{ij}\bigl(d(x_i,\sigma x_j)\bigr)$ and 
$b_{ij}= K_{ij}\bigl(d(x_i,x_j)\bigr) - K_{ij}\bigl(d(x_i,\sigma x_j)\bigr)$,
so that 
$$
K_{ij}\bigl(d(\sigma^{\eps_i}x_i,\sigma^{\eps_j}x_j)\bigr) 
= a_{ij} + b_{ij}\mathbf{1}_{\eps_i=\eps_j}\ .
$$
The product term in Eq.~(\ref{eq:Riesz-2pt-proof-1}) expands to
$$
\prod_{i<j} K_{ij}\bigl(d(\sigma^{\eps_i}x_i,\sigma^{\eps_j}x_j)\bigr)
= \sum_{\Gamma} \Bigl(\prod_{ij\not\in E} a_{ij}\Bigr)
\Bigl( \prod_{ij\in E} b_{ij} {\mathbf 1}_{\eps_i=\eps_j}\Bigr)
=: C_\Gamma {\mathbf 1}_{\eps_i=\eps_j, ij\in E}
\ ,
$$
where $\Gamma$ runs over all proper graphs on the 
vertex set $V=\{1,\dots, m\}$, and $E$ is the set of edges
of~$\Gamma$.  Inserting the expansion into 
Eq.~(\ref{eq:Riesz-2pt-proof-1}) and exchanging the order of
summation shows that each graph contributes a nonnegative term 
\begin{equation}
\label{eq:Riesz-2pt-proof-2}
C_\Gamma \sum_{\eps_i\in\{0,1\}, i\in V} 
F\bigl(u_1(\sigma^{\eps_1}x_1),\dots, u_m(\sigma^{\eps_m}x_m)\bigr)
{\mathbf 1}_{\eps_i=\eps_j, ij\in E}
\end{equation}
to the integral in Eq.~(\ref{eq:Riesz-2pt-proof-1}).
If $\Gamma$ is connected, then 
\begin{eqnarray}
\nonumber 
&&\hskip -3cm
\sum_{\eps_i\in\{0,1\}, i\in V} 
F\bigl(u_1(\sigma^{\eps_1}x_1),\dots, u_m(\sigma^{\eps_m}x_m)\bigr)
{\mathbf 1}_{\eps_i=\eps_j, ij\in E} \\
\nonumber 
&=& F\bigl(u_1(x_1),\dots, u_m(x_m)\bigr) 
+ F\bigl(u_1(\sigma x_1),\dots, u_m(\sigma x_m)\bigr) \\[0.1cm] 
\label{eq:Riesz-2pt-proof}
&\le& F\bigl(u_1^\sigma (x_1)),\dots, u_m^\sigma (x_m)\bigr) 
+  F\bigl(u_1^\sigma (\sigma x_1)),\dots, u_m^\sigma (\sigma x_m)\bigr)
\\[0.1cm] 
\nonumber &=& \sum_{\eps_i\in\{0,1\}, i\in V} 
F\bigl(u_1^\sigma (\sigma^{\eps_1}x_1),\dots, u_m^\sigma 
         (\sigma^{\eps_m}x_m)\bigr)
{\mathbf 1}_{\eps_i=\eps_j, ij\in E} \ ,
\end{eqnarray}
where the second step follows from Lemma~\ref{lem:Lorentz}
with $z_i=u_i(x_i)$ and $w_i=u_i(\sigma x_i)$.

If $\Gamma$ is not connected, choose a connected
component $\Gamma'$ and let $\Gamma''$ be its complement.
Let $E'$, $E''$, $V'$, and $V''$ be the corresponding edge and vertex sets.
The sum in Eq.~(\ref{eq:Riesz-2pt-proof-2}) can be decomposed~as
$$
\sum_{\eps_i\in \{0,1\}, i\in V''} 
\left\{ \sum_{\eps_i\in \{0,1\}, i\in V'} 
F\bigl(u_1(\sigma^{\eps_1}x_1),\dots, u_m(\sigma^{\eps_m}x_m)\bigr)
\mathbf{1}_{\eps_i=\eps_j, ij\in E'}\right\}
\mathbf{1}_{\eps_i=\eps_j, ij\in E''}\ .
$$
The key observation is that Eq.~(\ref{eq:Riesz-2pt-proof})
applies to the term in braces for fixed $\eps_i, i\in V''$; 
in other words, the contribution of $\Gamma$ can only increase if 
$u_i$ is replaced by $u_i^\sigma$ for all $i\in V'$.
An induction over the connected components of $\Gamma$ shows that
$$
\begin{array}{l}
\hskip -1cm  
\displaystyle{
\sum_{\eps_i\in\{0,1\}, i\in V} 
F\bigl(u_1(\sigma^{\eps_1}x_1),\dots, u_m(\sigma^{\eps_m}x_m)\bigr)
{\mathbf 1}_{\eps_i=\eps_j, ij\in E}}\\
\displaystyle{\hskip 2cm \le
\sum_{\eps_i\in\{0,1\}, i\in V} 
F\bigl(u_1^\sigma (\sigma^{\eps_1}x_1),\dots, u_m^\sigma 
                    (\sigma^{\eps_m}x_m)\bigr)
{\mathbf 1}_{\eps_i=\eps_j, ij\in E}}
\end{array}
$$
for any graph $\Gamma=(E,V)$. Adding  the contributions  of all graphs
shows that the integrand in Eq.~(\ref{eq:Riesz-2pt-proof-1})
increases pointwise under two-point rearrangement, 
and Eq.~(\ref{eq:Riesz-2pt}) follows.

\medskip{\em Equality statement:}  
Let $\Gamma_0$ be the graph defined in
the statement of the lemma, and let $E_0$ be its edge set.
By assumption,
$$
C_{\Gamma_0}
= \Bigl(\prod_{ij\not\in in E_0} 
K_{ij}\bigl(d(x_i,x_j)\bigr) - K_{ij}\bigl(d(\sigma x_i,x_j)\bigr) \Bigr)
\Bigl( \prod_{ij\in E_0} K_{ij}\bigl(d(\sigma x_i,x_j)\bigr)\Bigr)
>0\
$$
for $a.e.\ x_1,\dots, x_m\in H_+$.
If $\Delta_{ij}F>0$, then Lemma~\ref{lem:Lorentz} implies that 
Eq.~(\ref{eq:Riesz-2pt-proof}) is strict unless 
$$
\bigl(u_i(x_i)-u_i(\sigma x_i)\bigr)\bigl(u_j(x_j)-u_j(\sigma x_j)\bigr)
\ge 0\ , \quad a.e.\ x_i,x_j\in H_+\ .
$$
If $u_i$ and $u_j$ are not 
symmetric under $\sigma$, the product is not identically zero.
Since $x_i$ and $x_j$ can vary independently, 
this means that $u_i(x)-u_i(\sigma x)$ and
$u_j(x)-u_j(\sigma x)$ cannot change sign on $H_+$. 
We conclude that equality in Eq.~(\ref{eq:Riesz-2pt}) 
implies that either $u_i=u_i^\sigma$ and $u_j=u_j^\sigma$,
or $u_i=u_i^\sigma\circ\sigma$ and $u_j=u_j^\sigma\circ\sigma$.
\hfill $\Box$\quad

\bigskip\noindent 
Draghici also used Baernstein-Taylor approximation
to obtain Eq.~(\ref{eq:Riesz-ext}) from 
Eq.~(\ref{eq:Riesz-2pt}). If $F$ is bounded 
and continuous and $K_{ij}$ is bounded for $1\le i<j\le m$, 
then for bounded functions $u_1,\dots, u_m$ that are supported
in a common ball $B$ the inequality follows
from Lemma~\ref{lem:Riesz-2pt} by approximating the symmetric 
decreasing rearrangement with a sequence of two-point 
rearrangements, see Eq.~(\ref{eq:2pt-converge}).
Dominated convergence applies as in Eq.~(\ref{eq:dominated}),
since the integrations extend only over the bounded set $B^m$.

\section{Proof of the main results}
\label{sec:proofs}
\setcounter{equation}{0}


\noindent{\sc Proof of Theorem~\ref{thm:HL}.} \ 
{\em The inequality for Borel integrands:}
Let $F$ be a supermodular Borel function with
$F({\bf 0})=0$, and let 
and $u_1,\dots, u_m$ be nonnegative measurable functions that vanish
at infinity, as in the statement of the theorem. Denote by
$$
\I(u_1,\dots, u_m) := \int_\XX F(u_1(x),\dots, u_m(x))\, dx
$$
the left hand side of Eq.~(\ref{eq:HL-ext}). 
Replacing $F({\bf y})$ by $F({\bf y})  - \sum_{i=1}^m F(y_i{\bf e}_i)$
and using that $F(u_i(\cdot)\,{\bf e}_i)$ and
$F(u_i^*(\cdot)\,{\bf e}_i)$ contribute equally to the two sides
of Eq.~(\ref{eq:HL-ext}), we may assume $F$ to be
nondecreasing in each variable.

Fix $L>0$, and replace $u_i$ by the bounded function
$$
u_i^L(x) := \min\left\{u_i(x),L\right\} {\mathbf 1}_{\{|x|<L\}}  
$$
for $i=1,\dots, m$.  Then
\begin{equation} \label{eq:HL-proof-1}
F(u_1^L, \dots,u_m^L)= F^L(u_1^L, \dots,u_m^L)\ ,
\end{equation}
where $F^L$ is the function defined in Lemma~\ref{lem:cutoff}. 
By construction, $F^L$ is bounded, and by Lemma~\ref{lem:cutoff}
it is nondecreasing and supermodular.
By Lemma~\ref{lem:Sklar},
there exist nondecreasing functions $\phi_i$ with $\phi_i(0)=0$
and a continuous supermodular function
$\tilde F^L$ on $\RR^m_+$ such that
\begin{equation} \label{eq:factorize}
F^L(y_1,\dots, y_m)=\tilde F^L(\phi_1(y_1),\dots, \phi_m(y_m))\ .
\end{equation}
Since $\phi_i$ is nondecreasing and vanishes at zero, $u_i^L$ 
is compactly supported, and $(u_i^L)^*\le (u_i^*)^L$ pointwise 
by construction, we have 
\begin{equation} \label{eq:phi-*}
(\phi_i\circ u_i^L)^* = \phi_i\circ (u_i^L)^* \le 
       \phi_i\circ (u_i^*)^L
\end{equation}
for $i=1,\dots,m$. By Theorem 1 of~\cite{Brock}),
$$
\int_\XX \tilde F^L(\phi_1\circ u_1^L(x),\dots, \phi_m\circ u_m^L(x))\, dx
\le
\int_\XX \tilde F^L((\phi_1\circ u_1^L)^*(x),\dots, 
(\phi_m\circ u_m^L)^*(x))\, dx\ .
$$
With Eqs.~(\ref{eq:HL-proof-1})-(\ref{eq:phi-*}), this becomes
$$
\I(u_1^L,\dots, u_m^L)\le \I\bigl((u_1^*)^L,\dots, (u_m^*)^L\bigr)\ .
$$
Since $u_i^L(x) = u_i(x)$ for $L\ge \max\{u_i(x), |x|\}$,
we see that $F(u_i^L(x),\dots, u_m^L(x))$ converges
pointwise to $F(u_1(x),\dots, u_m(x))$,
and Eq.~(\ref{eq:HL-ext}) follows by monotone convergence.

\medskip {\em Equality statement:} \
Combining 
Eq.~(\ref{eq:HL-2pt}) with Eq.~(\ref{eq:HL-ext})
and using that $u_i^\sigma$ is equimeasurable with $u_i$, 
we see that
$$
\I(u_1,\dots, u_m)\le  \I(u_1^\sigma,\dots, u_m^\sigma) \le 
\I(u_1^*,\dots, u_m^*)\ .
$$
Hence equality in Eq.~(\ref{eq:HL-ext}) implies 
equality in Eq.~(\ref{eq:HL-2pt}) 
for every choice of the reflection $\sigma$. Given two points
$x,x'$ in $\XX$, choose $\sigma$ such that $\sigma(x)=x'$.  
If $\Delta_{ij}F>0$ for some $i\ne j$, then
$u_i(x)-u_i(x')$ and $u_j(x)-u_j(x')$ cannot have opposite
signs by Lemma~\ref{lem:HL-2pt}. 
If $u_i=u_i^*$ is strictly radially decreasing, then
it follows that $u_j^\sigma=u_j$ for every reflection  $\sigma$
that does not fix $x^*$.  By Eq.~(\ref{eq:2pt-identify}), 
$u_j=u_j^*$ as claimed.
\hfill$\Box$\quad

\bigskip\noindent{\sc Proof of Theorem~\ref{thm:Riesz}.} \ 
{\em The inequality for Borel integrands:} 
The proof of Eq.~(\ref{eq:Riesz-ext}) proceeds along the same lines
as the proof of Eq.~(\ref{eq:HL-ext}). Let
$$
\I(u_1,\dots, u_m) := 
\int_\XX\cdots \int_\XX F(u_1(x_1),\dots , u_m(x_m))\, 
   \prod_{i<j} K_{ij}(d(x_i,x_j))\, dx_1\dots dx_m 
$$
be the left hand side of Eq.~(\ref{eq:Riesz-ext}).
As before, we may assume that $F$ is nondecreasing in each variable. 
We replace $F$ with $\tilde F^L$, $u_i$ with $\phi_i\circ u_i^L$, 
$K_{ij}$ with $K_{ij}^L=\min\{K_{ij},L\}$, and set
$$
\I^L(u_1,\dots, u_m) := \int_{\XX}\cdots \int_{\XX} F^L\bigl(u_1(x_1),\dots , 
u_m(x_m)\bigr)\, 
  \prod_{i<j} K_{ij}^L\bigl(d(x_i,x_j)\bigr)\, dx_1\dots dx_m\ .
$$
Applying Theorem 2.2 of~\cite{Draghici}, we obtain with
the help of
Eqs.~(\ref{eq:HL-proof-1})-(\ref{eq:phi-*})
$$
\I^L(u_1^L,\dots,u_m^L) \le \I^L\bigl((u_1^*)^L,\dots,(u_m^*)^L\bigr)\ .
$$
Eq.~(\ref{eq:Riesz-ext}) follows by taking 
$L\to\infty$ and using monotone convergence.

\medskip {\em Equality statement:}
Consider the set $S_i$ of all reflections $\sigma$ of 
$\XX$ that fix $u_i$.  
If $u_i$ is non-constant, then $S_i$ is a closed proper subset 
of the space of all reflections on $\XX$. This subset is nowhere dense, 
since any open set of reflections generates the entire isometry
group of $\XX$.  If Eq.~(\ref{eq:Riesz-ext}) 
holds with equality, then the two-point rearrangement inequality in 
Eq.~(\ref{eq:Riesz-2pt}) holds with equality for every 
reflection $\sigma$. For $\sigma\not \in  S_i$, 
Lemma~\ref{lem:Riesz-2pt} implies that either $u_j=u_j^\sigma$ 
or $u_j=u_j^\sigma\circ\sigma$.
Since $S_i$ is nowhere dense, it follows from the continuous dependence of
$u^\sigma$ on $\sigma$ that $u_j$ agrees with
either $u_j^\sigma$ or $u_j^\sigma\circ\sigma$
also for $\sigma\in S_i$. By Eq.~(\ref{eq:2pt-identify-tau}),
there exists a translation $\tau$ such that $u_j=u_j^*\circ\tau$.
Lemma~\ref{lem:Riesz-2pt} implies furthermore that
$u_i$ agrees with $u_i^\sigma$ when
$u_j=u_j^\sigma$, and with 
$u_i^\sigma\circ\sigma$ when $u_j=u_j^\sigma\circ\sigma$. 
We conclude that $u_i=u_i^*\circ\tau$.
\hfill$\Box$\quad

\section{Concluding remarks}
\label{sec:conc}
\setcounter{equation}{0}

In the proof of Eq.~(\ref{eq:Riesz-ext}) and its two-point version
in Eq.~(\ref{eq:Riesz-2pt}), the kernels $K_{ij}$ 
played a very different role from the functions $u_1,\dots,u_m$
that enter into the integrand. However, the Riesz functional 
on the left hand side of Eq.~(\ref{eq:Riesz}) depends
equally on  $u$, $v$, and $w$.  We will use the connection of Riesz' 
inequality with the Brunn-Minkowski inequality to construct examples 
where the two-point rearrangement fails for Eq.~(\ref{eq:Riesz}).

The Brunn-Minkowski inequality says that the measures 
of two subsets $A,B\subset\RR^n$ are related to the measure 
of their Minkowski sum $A+B= \{a+b:\ a\in A, b\in B\}$ by
$$\lambda(A)^{1/n} + \lambda(B)^{1/n} \le \lambda(A+B)^{1/n}\ . $$
Recognizing the two sides of the inequality as proportional
to the radii of the balls $A^*+B^*$ and $(A+B)^*$,
we rewrite it as the rearrangement inequality
\begin{equation}
\label{eq:BM}
\lambda (A^*+B^*) \le \lambda (A+B)\ .
\end{equation}
Eq.~(\ref{eq:BM}) follows rather directly from Riesz' 
inequality in Eq.~(\ref{eq:Riesz}), because
the support of the convolution of two 
nonnegative functions is essentially the Minkowski 
sum of their supports.  Conversely, the
Brunn-Minkowski inequality enters into the proof
of the Brascamp-Lieb-Luttinger inequality~\cite{BLL},
of which Eqs.~(\ref{eq:Riesz}) and~(\ref{eq:Riesz-ext}) 
are special cases.

Equality in the Brunn-Minkowski inequality
implies that $A$ and $B$ differ only by sets of measure zero from 
two independently scaled and translated 
copies of a convex body~\cite{H-Oh}. 
Let $A=B$ be an ellipsoid in $\RR^n$ with $n>1$ that 
is centered at a point $c\ne 0$, so that
Eq.~(\ref{eq:BM}) holds with equality.
If $\sigma$ is the reflection at a hyperplane through 
$c$ that is not a hyperplane of symmetry
for $A$ and $B$, then $A^\sigma$ and $B^\sigma$ are
non-convex, and  therefore
$$
\lambda(A^\sigma+B^\sigma) > \lambda(A^*+B^*)= \lambda(A+B)\ .
$$
Choosing $u$, $v$, and $w$ as the characteristic functions of $A$, 
$A+B$, and~$B$ provides an example where the Riesz functional strictly 
decreases under two-point rearrangement. For an example 
of this phenomenon in one dimension, 
consider the symmetric decreasing functions 
$$
u(x)={\mathbf 1}_{|x-2|<\eps}\ ,\quad v(x)=w(x)={\mathbf 1}_{|x-1|<\eps}\ ,
$$
and let $\sigma$ be the reflection at $x=1$. Then
$$
u^\sigma(x)={\mathbf 1}_{|x|<\eps}\ ,
\quad v^\sigma(x)=w^\sigma(x)={\mathbf 1}_{|x-1|<\eps}\ ,
$$
and if $0<\eps\le \frac{1}{2}$, Riesz' inequality fails for $\sigma$,
$$
\begin{array}{lcl}
\displaystyle{\int_\RR\int_\RR u(x)v(x')w(x-x')\, dx dx' }
&>& 0 \\
&=& 
\displaystyle{\int_\RR\int_\RR u^\sigma(x)v^\sigma(x)w^\sigma(x-x')\, dxdx'\ .}
\end{array}
$$

While the two-point rearrangement is not useful
for Eq.~(\ref{eq:Riesz}), the layer-cake representation
of Crowe-Zweibel-Rosenbloom shows that 
\begin{equation}
\label{eq:Riesz-ext-3}
\begin{array}{l}
\displaystyle{\hskip -2cm
\int_{\RR^n}\int_{\RR^n} F\bigl(u(x),v(x'),w(x-x')\bigr)\, dxdx'
}\\
\displaystyle{ \qquad 
\le \ \int_{\RR^n}\int_{\RR^n} F\bigl(u^*(x),v^*(x'),w^*(x-x')\bigr)\, dxdy}
\end{array}\end{equation}
for any integrand that can be written as the joint distribution function
of a Borel measure $\mu_F$ on $\RR^3_+$,
$$
F(y_1,y_2,y_3)= \mu_F\bigl( [0,y_1)\times [0,y_2)\times [0,y_3)\bigr)\ .
$$
Such integrands are left continuous,
vanish at the origin, and satisfy $\Delta_{i_1,\dots,i_\ell}F\ge 0$ for
every choice of $\ell\le 3$ distinct indices. 
Lemma~\ref{lem:Sklar} allows 
to accommodate integrands in Eq.~(\ref{eq:Riesz-ext-3})
that are only Borel measurable. 
The main condition is that $\Delta_{123}F\ge0$;
the second-order monotonicity conditions can be replaced by integrability 
assumptions on the negative part $F_-$ similar to 
Eq.~(\ref{eq:Riesz-integrable}). 
To ensure that the functional is finite at least when $u,v,w$ are bounded
and compactly supported, $F$ should vanish on the coordinate axes.
For example, Eq.~(\ref{eq:Riesz-ext-3}) holds for
$$F(u,v,w)= \frac{uvw}{(1+u)(1+v)(1+w)} - (uv+uw+vw) $$
since $\Delta_{123}F> 0$, even though $\Delta_{ij}F< 0$ 
for all $i\ne j$. 

For Borel integrands satisfying
$\Delta_{123}F>0$, equality in Eq.~(\ref{eq:Riesz-ext-3})
implies that every triple of level sets of $u,v,w$ produces 
equality in Eq.~(\ref{eq:Riesz}).  These equality cases were
described in~\cite{Burchard}.  In particular, if two of the three 
functions $u,v,w$ are known to have continuous distribution
functions and the value of the functional is finite, 
then equality implies that $u,v,w$ are equivalent 
to $u^*, v^*, w^*$ under the symmetries of the functional
(see Theorem~2 of~\cite{Burchard}).

By the same line of reasoning, the Brascamp-Lieb-Luttinger 
inequality~\cite{BLL} implies that
$$
\I(u_1,\dots,u_m) :=\int_{\RR^n}\dots\int_{\RR^n}
F\Bigl(u_1\bigl(\sum_{j=1}^k a_{1j}x_j\bigr),\dots,
u_m\bigl(\sum_{j=1}^k a_{mj}x_j\bigr)\Bigr)\,
dx_1 \dots dx_k
$$
increases under symmetric decreasing rearrangement, if 
$\Delta_{i_1,\dots,i_\ell}F\ge 0$ for all choices of distinct 
indices $i_1,\dots,i_\ell$ with $\ell\le m$. 
Interesting examples are integrands of the form 
in Eq.~(\ref{eq:Morpurgo}), where $\Phi$ is completely monotone 
in the sense that all its distributional derivatives are nonnegative. 
If $\Delta_{i_1\dots i_\ell}F>0$ for all choices of 
$i_1,\dots,i_\ell$, then the last statement of 
Lemma~\ref{lem:cutoff} can be used to show that 
the extended Brascamp-Lieb-Luttinger inequality
has the same equality cases  as the original inequality.
The characterization of these equality cases remains an open problem.



\end{document}